\documentclass[preprint,12pt]{elsarticle}

\usepackage{natbib}
\usepackage{amssymb}
\usepackage{color}
\newtheorem{theorem}{Theorem}[section]
\newtheorem{definition}{Definition}[section]
\newtheorem{proposition}{Proposition}[section]
\newtheorem{lemma}{Lemma}[section]

\newtheorem{remark}{Remark}[section]

\journal{IMA J. Math. Control Info.}

\begin{document}

\begin{frontmatter}

\title{Regional boundary controllability of time fractional  diffusion processes}

\author[dhu]{Fudong Ge}
\ead{gefd2011@gmail.com}
\author[ucmerced]{YangQuan Chen\corref{cor}}
\ead{ychen53@ucmerced.edu}
\author[dhutwo]{Chunhai Kou}
\ead{kouchunhai@dhu.edu.cn}
\cortext[cor]{Corresponding author}

\address[dhu]{College of Information Science and Technology, Donghua University, Shanghai 201620, PR China}
\address[ucmerced]{Mechatronics, Embedded Systems and Automation Lab, University of California, Merced, CA 95343, USA}
\address[dhutwo]{Department of Applied Mathematics, Donghua University, Shanghai 201620, PR China}

\begin{abstract}
In this paper, we are concerned with the regional boundary controllability of the Riemann-Liouville time fractional
diffusion systems of order $\alpha\in (0,1]$. The characterizations of strategic
actuators are established when the systems studied are regionally boundary controllable. The determination of
control to achieve regional boundary controllability with minimum energy is explored.  We also show a connection
between the regional internal controllability and regional boundary controllability. Several useful results for
the optimal control from an implementation  point of view are presented  in the end.
\end{abstract}
\par
\end{frontmatter}

\section{Introduction}
\setcounter{section}{1}\setcounter{equation}{0}
In the past several decades, a lot of work has been carried out to deal with the problem of steering a system to a target state,
especially after the introduction of the notions of actuators and sensors \cite{AEIJAI,AEIJAI2}. However, in many real-world
applications, we are only concerned with those cases where the target states of the problem studied are defined in a given
subregion of the whole space domain. Then the regional idea emerges and we refer the reader to \cite{217,zerrik,179} for more
information on the concept of regional analysis for the Gaussian diffusion process.  Besides, it should be
pointed out that not only does the concept of regional analysis make sense closer to real-world problems, it also
generalizes the results of existence contributions.

In addition, after the introduction of continuous time random walks (CTRWs) by Montroll
and Weiss \cite{CTRWfirst}, the anomalous diffusion equation of fractional order has attracted increasing interest and
has been proven to be a useful tool in modeling many real-world problems \cite{ralf,ralf2,Cartea,3M,Adams}.
More precisely, the mean squared displacement (MSD) of
anomalous diffusion  process is described by a power law of fractional exponent, which is  smaller (in the case of
sub-diffusion) or bigger (in the case of super-diffusion) than that of Brownian motion. It is confirmed that
the time fractional  diffusion system, where the traditional first order time derivative is replaced by a
Riemann-Liouville time fractional derivative of order $\alpha\in (0,1]$, can be used to well characterize those
sub-diffusion process \cite{ralf,ralf2}.  For example, the flow through porous media microscopic processes
\cite{Uchaikin}, or
swarm of robots moving  through dense forest \cite{2Spears} etc.
For the fractional calculus, as we all know, it has shown great potential
in science and engineering applications and some phenomena such as self-similarity, nonstationary,
non-Gaussian process and short or long memory process are all closely related
to fractional calculus \cite{opAgrawal,Podlubny,Kilbas}. It is now widely believed that, using fractional calculus in
modeling can better capture the complex dynamics of natural and man-made systems, and
fractional order controls can offer better performance not achievable before using integer
order controls \cite{Torvik,Mandelbrot}, which in fact raise important theoretical challenges and open new research opportunities.

Motivated by the argument above, the contribution of this present work is on the regional boundary controllability of the anomalous transport process described by
time fractional diffusion systems. More precisely, for an open bounded subset $\Omega\subseteq \mathbf{R}^n$ with smooth boundary
$\partial\Omega$, we consider:
\begin{itemize}
 \item{} A subregion $\Gamma$ of $\partial\Omega$ which may be unconnected.
\vspace{-0.35cm}
 \item{}  Various kinds of actuators (zone, pointwise, internal or boundary) acting as controls.
\end{itemize}
The rest of this paper is organized as follows. The mathematical concept of regional boundary
controllability and several preliminaries are presented in the next section, then we present an
example which is regional boundary controllability but not globally boundary controllable. Section $\ref{rbcsection3}$ is focused
on the characterizations of $\Gamma-$strategic actuators and our main result on regional boundary
controllability with minimum energy problem is given in Section $\ref{rbcsection4}$. In Section $\ref{rbcsection5}$, a connection between internal and
boundary regional controllability is established and at last, we work out some useful results for
the optimal control from an implementation point of view.

\section{Regional boundary controllability}
\label{rbcsection2}
\setcounter{section}{2}\setcounter{equation}{0}
\subsection{Problem statement}
In this paper,  we consider the following abstract time fractional diffusion system:
\begin{equation}\label{problem}
\left\{\begin{array}{l}
_0D^{\alpha}_tz(t)=Az(t)+Bu(t), ~~ t\in [0,b],
\\ \lim\limits_{t\to 0^+} {}_0I^{1-\alpha}_t z(t)=z_0,
\end{array}\right.
\end{equation}
where $A$ generates a strongly continuous semigroup $\{\Phi(t)\}_{t\geq 0}$ on the Hilbert
space $Z:=H^1(\Omega)$, $z\in L^2(0,b;Z)$  and the initial vector $z_0\in Z$. It is supposed
that  $B: \mathbf{R}^p\to Z$ is the control operator and $u\in L^2(0,b;\mathbf{R}^p)$ depends
on the number and structure of actuators.
Moreover, the
Riemann-Liouville fractional derivative $_0D_t^\alpha$ and the Riemann-Liouville fractional
integral $_0I^{\alpha}_t$  are, respectively, given by \cite{Podlubny},\cite{Kilbas}
\begin{equation}\label{caputoderivative}
_0D_t^\alpha z(t)=\frac{d}{d t}{}_0I_t^{1-\alpha}z(t),~~\alpha\in (0,1]~\mbox{ and }~{}_0I^{\alpha}_tz(t)=\frac{1}{\Gamma(\alpha)}\int^{t}_{0}{(t-s)^{\alpha-1}z(s)ds},~~\alpha>0.
\end{equation}

\begin{definition} \cite{GeFDTA-rc}\label{lemma2.1} For any given $f\in  L^2\left(0,b; Z\right),$ $\alpha\in (0,1],$
a function $v\in L^2\left(0,b; Z\right)$  is said to be a mild solution of the following system
\begin{equation}\label{problem2}
\left\{\begin{array}{l}
{}_0D^{\alpha}_tv(t)=Av(t)+f(t),~t\in [0,b],
\\ \lim\limits_{t\to 0^+}{}_0I_t^{1-\alpha}v(t)=v_0\in Z,
\end{array}\right.
\end{equation}
if it satisfies
\begin{eqnarray}
z(t)=K_\alpha(t)v_0+\int_0^t(t-s)^{\alpha-1}{K_\alpha(t-s)f(s)}ds,
\end{eqnarray}
where
$
K_\alpha(t)=\alpha\int_0^\infty{\theta\phi_\alpha(\theta)\Phi(t^\alpha\theta)}d\theta$, $\{\Phi(t)\}_{t\geq 0}$ is the strongly continuous semigroup generated by $A$, $\phi_\alpha(\theta)=
\frac{1}{\alpha}\theta^{-1-\frac{1}{\alpha}}\psi_\alpha(\theta^{-\frac{1}{\alpha}})$ and  $\psi_\alpha$ is a probability density function defined by
\begin{eqnarray}\label{pdf}
\psi_\alpha(\theta)=\frac{1}{\pi}\sum\limits_{n = 1}^\infty {( - 1)^{n - 1}\theta^{-\alpha n-1}\frac{\Gamma(n\alpha+1)}{n!}
\sin(n\pi \alpha)}, {\kern 6pt}\theta >0. \end{eqnarray}
In addition, we have \cite{Mainardi,zhouyong}
\begin{eqnarray}\int_0^\infty{\psi_\alpha(\theta)}d\theta=1 ~\mbox{ and }~ \int_0^\infty{\theta^\nu
\phi_\alpha(\theta)}d\theta=\frac{\Gamma(1+\nu)}{\Gamma(1+\alpha\nu)}, {\kern 4pt}\nu\geq 0.\end{eqnarray}
\end{definition}
By Lemma $\ref{lemma2.1}$, the mild solution $z(.,u)$ of $(\ref{problem})$ can be  given by
\begin{eqnarray}
z(t,u)=K_\alpha(t)z_0+\int_0^t(t-s)^{\alpha-1}{K_\alpha(t-s)Bu(s)}ds.
\end{eqnarray}
Let $H:L^2(0,b;\mathbf{R}^p) \to  Z$  be
 \begin{eqnarray}Hu=\int_0^b(b-s)^{\alpha-1}{K_\alpha(b-s)Bu(s)}ds,~\forall u\in L^2(0,b;\mathbf{R}^p)_.\end{eqnarray}
Suppose that $\{\Phi^*(t)\}_{t\geq 0}$, generated by the adjoint operator of $A$, is also a strongly continuous
semigroup on the space $Z$.
For any $v\in Z,$ it follows from  $\left<Hu,v\right>=\left<u,H^*v\right>$ that
\begin{eqnarray}\label{H*}
H^* v=B^*(b-s)^{\alpha-1}K_\alpha^*(b-s)v,
\end{eqnarray} where $\left<\cdot,\cdot\right>$ is
the duality pairing of space $Z$, $B^*$ is the adjoint operator of $B$ and
\[K^*_\alpha(t)= \alpha\int_0^\infty{\theta\phi_\alpha(\theta)\Phi^*(t^\alpha\theta)}d\theta.\]

Let $\gamma:H^1(\Omega)\to H^{\frac{1}{2}}(\partial\Omega)$ be the trace operator of order zero, which is linear continuous
and surjective, $\gamma^*$ denotes the adjoint operator. Moreover, if $\Gamma \subseteq \partial \Omega,$
$p_\Gamma:H^{\frac{1}{2}}(\partial \Omega)\to H^{\frac{1}{2}}(\Gamma)$ defined by
\begin{eqnarray}
p_\Gamma z:=z|_\Gamma
\end{eqnarray}
and for any $\bar{z}\in H^{\frac{1}{2}}(\Gamma),$
 the adjoint operator $p_\Gamma^*$ can be given by
\begin{eqnarray}
p_\Gamma^* \bar{z}(x):=\left\{\begin{array}{l}\bar{z}(x), ~~x\in \Gamma ,  \\0,~~~x\in \partial\Omega\backslash \Gamma.\end{array}\right.
\end{eqnarray}


\subsection{Definition and characterizations}

Let $\omega \subseteq \Omega $ be a given region of positive Lebesgue measure. Denote the projection operator on $\omega$ by the restriction map
\begin{eqnarray}p_\omega: H^1(\Omega)\to H^1(\omega),\end{eqnarray}
then we are ready to state the following definitions.

\begin{definition} \label{rcdefinition1}
 The system $(\ref{problem})$ is said to be exactly (respectively, approximately) regionally controllable on $\omega$ at time $b$
 if for any $y_b\in H^1(\omega)$, given $\varepsilon>0,$ there exists a control $u\in L^2(0,b; \mathbf{R}^p)$  such that
\begin{eqnarray}
p_\omega z(b,u)=y_b\left(\mbox{respectively,  }\|p_\omega z(b,u)-y_b \|_{H^1(\omega)}\leq \varepsilon\right).
\end{eqnarray}
\end{definition}

 \begin{definition} \label{rbcdefinition1}
The system $(\ref{problem})$ is said to be exactly (respectively, approximately) regionally boundary
controllable on $\Gamma\subseteq \partial \Omega$ at time $b$ if for any $z_b\in H^{\frac{1}{2}}(\Gamma)$, given $\varepsilon>0,$ there exists
 a control $u\in L^2(0,b; \mathbf{R}^p)$  such that
\begin{eqnarray}
p_\Gamma \left(\gamma z(b,u)\right)=z_b\left(\mbox{respectively,  }\|p_\Gamma \left(\gamma z(b,u)\right)-z_b \|_
{H^{\frac{1}{2}}(\Gamma)}\leq \varepsilon\right).
\end{eqnarray}
\end{definition}

\begin{proposition}\label{proposition1} The following properties are equivalent:

$(1)$ The system $(\ref{problem})$ is exactly regionally boundary controllable on $\Gamma$  at time $b$;

$(2)$ $Im( p_\Gamma \gamma  H) =H^{\frac{1}{2}}(\Gamma)$;

$(3)$ $ Ker(p_\Gamma )+ Im( \gamma  H)= H^{\frac{1}{2}}(\partial \Omega);$

$(4)$ For any $z\in  H^{\frac{1}{2}}(\Gamma),$ there exists a positive constant $c$ such that
\begin{eqnarray}
\|z\|_{ H^{\frac{1}{2}}(\Gamma)}\leq c\|H^*\gamma^*p_\Gamma^* z\|_{L^2(0,b;\mathbf{R}^p)}.
\end{eqnarray}
\end{proposition}
\textbf{Proof.}
By Definition $\ref{rbcdefinition1}$, it is not difficult to see that $(1)\Leftrightarrow (2).$

$(2)\Rightarrow (3):$ For any $z\in H^{\frac{1}{2}}(\Gamma),$ let $\hat{z}$ be the extension of $z$ to $ H^{\frac{1}{2}}(\partial\Omega).$
Since  $Im( p_\Gamma \gamma  H )= H^{\frac{1}{2}}(\Gamma)$, there exists $u\in L^2(0,b;\mathbf{R}^p)$, $z_1\in Ker( p_\Gamma)$ such that $\hat{z}=z_1+\gamma Hu.$

$(3)\Rightarrow (2):$ For any $\tilde{z}\in H^{\frac{1}{2}}(\partial \Omega)$, $\tilde{z}=z_1+z_2,$ where $z_1\in Ker( p_\Gamma)$ and $z_2\in  Im(\gamma H)$.
Then there exists a control $u\in L^2(0,b;\mathbf{R}^p)$ such that $\gamma Hu=z_2$. Hence, it follows from the definition of $p_\Gamma$ that $Im( p_\Gamma\gamma  H)=H^{\frac{1}{2}}(\Gamma)$.

$(1)\Leftrightarrow (4):$ The equivalence between $(1)$ and $(4)$ can be deduced from the following general result \cite{dualityrelationship}:
Let $E, F, G $ be reflexive Hilbert spaces and $f \in \mathcal{L}(E,G),$ $g \in \mathcal{L}(F,G).$ Then
the following two  properties are equivalent

$(1)$  $Im(f )\subseteq Im( g),$

$(2)$ $\exists ~\gamma > 0$ such that $\|f^*z^*\|_{E^*}\leq \gamma\|g^*z^*\|_{F^*},~~\forall z^*\in G$.\\
By choosing $E=G= H^{\frac{1}{2}}(\Gamma),$ $F=L^2(0,b;\mathbf{R}^p),$ $f=Id_{ H^{\frac{1}{2}}(\Gamma)}$ and $g=p_\Gamma \gamma H$, then we complete the proof.

\begin{proposition} There is an equivalence among the following properties:

$\left<1\right>$ The system $(\ref{problem})$ is approximately regionally boundary controllable on $\Gamma$  at time $b$;

$\left<2\right>$ $\overline{Im( p_\Gamma \gamma  H )}=  H^{\frac{1}{2}}(\Gamma);$

$\left<3\right>$ $ Ker (p_\Gamma) + \overline{Im (\gamma H )}=  H^{\frac{1}{2}}(\partial\Omega);$

$\left<4\right>$ The operator $p_\Gamma \gamma HH^*\gamma ^*p^*_\Gamma$ is positive definite.
\end{proposition}
\textbf{Proof.}
By Proposition $\ref{proposition1}$,   $\left<1\right>\Leftrightarrow \left<2\right>\Leftrightarrow \left<3\right>$.
Finally, we show that $\left<2\right>\Leftrightarrow \left<4\right>$.
In fact, since
\begin{eqnarray*}
\overline{Im (p_\Gamma \gamma  H )}=  H^{\frac{1}{2}}(\Gamma)\Leftrightarrow (p_\Gamma \gamma  Hu,z)_{H^{1/2}(\Gamma)}=0 \mbox{ for any } u\in L^2(0,b;\mathbf{R}^p) \mbox{  implies  } z=0,
\end{eqnarray*}
where  $\left(\cdot,\cdot\right)_{H^{1/2}(\Gamma)}$ is the inner product of  $H^{\frac{1}{2}}(\Gamma).$  Let $u=H^*\gamma ^*p_\omega^* z$. Then we see that
\begin{eqnarray*}
\overline{Im (p_\Gamma \gamma  H) }=  H^{\frac{1}{2}}(\Gamma)\Leftrightarrow (p_\Gamma \gamma HH^*\gamma ^*p_\Gamma^* z,z)_{H^{1/2}(\Gamma)}=0 \mbox{  implies  } z=0,~z\in H^{\frac{1}{2}}(\Gamma),
\end{eqnarray*}
i.e.,  the operator $p_\Gamma \gamma HH^*\gamma ^*p^*_\Gamma$ is positive definite and the proof is complete.
\begin{remark}

$(1)$ A system which is boundary controllable on $\Gamma$ is   boundary
controllable on $\Gamma_1$ for every $\Gamma_1\subseteq \Gamma.$

$(2)$ The definitions $(\ref{rcdefinition1})$ can be applied to the case
where $\Gamma=\partial\Omega$ and  there exist systems that are not boundary controllable but which are regionally
 boundary controllable. This is illustrated by  the following example
\end{remark}

\subsection{An example}

Consider the following two dimension time fractional diffusion equation defined on $\Omega=[0,1]\times [0,1]$, which is excited by a zone actuator:
\begin{eqnarray}\label{examplel}
\left\{\begin{array}{l}
_0D^{\alpha}_{t}z(x,y,t)=\frac{\partial^2}{\partial x^2}z(x,y,t)+\frac{\partial^2}{\partial y^2}z(x,y,t)+p_D u(t)~ \mbox{ in } ~\Omega\times [0,b],
\\ \lim\limits_{t\to 0^+} {}_0I^{1-\alpha}_{t}z(x,y,t)=0~~ \mbox{ in }~ ~ \Omega,
\\ z(\xi,\eta, t)= 0~~ \mbox{ on } ~~\partial\Omega \times [0,b],
\end{array}\right.
\end{eqnarray}
where $\alpha\in (0,1],$
 $D=\{0\}\times[d_1,d_2]\subseteq \Omega$, $A=\frac{\partial^2}{\partial x^2}+\frac{\partial^2}{\partial y^2}$ with
$\lambda_{ij}=-(i^2+j^2)\pi^2 ,~\xi_{ij}(x,y)=2a_{ij}\cos(i\pi x)\cos(j\pi y)$, $a_{ij}=(1-\lambda_{ij})^{-\frac{1}{2}}$,
$
\Phi(t)z=\sum\limits_{i,j=1}^{\infty}{\exp(\lambda_{ij}t)(z,\xi_{ij})_{Z}\xi_{ij}}
$
and
$
K_\alpha(t)z(x)={\alpha}\int_0^\infty{\theta\phi_{\alpha}(\theta)\Phi(t^{\alpha}\theta)z(x)}d\theta
=\sum\limits_{i,j=1}^{\infty}E_{\alpha,\alpha}(\lambda_{ij}t^{\alpha})(z,\xi_{ij})_{Z}\xi_{ij}(x).
$
Further, since
\begin{eqnarray*}
(H^*\gamma^*z)(t)=(b-t)^{\alpha-1}\sum\limits_{i,j=1}^{\infty}E_{\alpha,\alpha}(\lambda_{ij}(b-t)^{\alpha})
(\gamma^*z,\xi_{ij})_{Z}(p_D,\xi_{ij})_{Z}
\end{eqnarray*}
and  $(p_D,\xi_{ij})_{Z}=\frac{2a_{ij}}{j\pi}\left[\sin(j\pi d_2)-\sin(j\pi d_1)+j\pi(\cos(j\pi d_2)-\cos(j\pi d_1)\right]$,  there exists $d_1,d_2\in [0,1]$ satisfying $ Ker (H^*) \neq
\{0\}$ $(\overline{Im (p_D H)}\neq L^2(\omega))$, i.e., the system $(\ref{examplel})$
is not boundary controllable.

Moreover, let $d_1=0$, $d_2=\frac{1}{2}$, $\Gamma=\{0\}\times[\frac{1}{4}, \frac{3}{4}]$ and
$z_*=\xi_{ij}(0,y), (i,j=4k,k=1,2,3,\cdots)$. Obviously, $z_*$ is not reachable on $\partial\Omega.$
However, since  \[E_{\alpha,\alpha}(t)>0~(t\geq 0)\mbox{ and } (p_D,\xi_{ij})_{Z}=\frac{2a_{ij}}{j\pi}\left[\sin(j\pi/2)+j\pi(\cos(j\pi/2)-1)\right],~j=1,2,\cdots,\]
we see that
\begin{eqnarray}\begin{array}{l}
(H^* \gamma^* p_\Gamma^*  z_*)(t)
=\sum\limits_{i,j=1}^\infty\frac{E_{\alpha,\alpha}(\lambda_{ij}(b-t)^{\alpha})}{(b-t)^{1-\alpha}}(\xi_{ij},\gamma^*z_*)_{H^{1/2}(\Gamma)}(p_D,\xi_{ij})_{Z}
\\{\kern 72pt}
=\sum\limits_{i,j=1,j\neq 4k}^\infty  \frac{2a_{ij}E_{\alpha,\alpha}(\lambda_{ij}(b-t)^{\alpha})}{j\pi(b-t)^{1-\alpha}}
(\xi_{ij},\gamma^*z_*)_{H^{1/2}(\Gamma)} \\ {\kern 102pt}\times \left[\sin(j\pi/2)+j\pi(\cos(j\pi/2)-1)\right]
\\ {\kern 72pt}\neq  0.
\end{array}\end{eqnarray}
Hence $z_*$ is regionally boundary  controllable on $\Gamma=\{0\}\times [\frac{1}{4}, \frac{3}{4}].$

To end this section, we finally recall a necessary lemma to be used afterwards.

\begin{lemma}\cite{Bernard}\label{infnitydifferential} Let $\Omega \subseteq \mathbf{R}^n$ be
 an open set and $C_0^{\infty}(\Omega)$ be the class of  infinitely differentiable functions on $\Omega$
 with compact support in $\Omega$ and $u\in L^1_{loc}(\Omega)$ be such that
\begin{eqnarray}
\int_\Omega u(x)\psi(x)dx=0,~~~\forall\psi \in C_0^{\infty}(\Omega).
\end{eqnarray}
Then $u=0$ almost everywhere in $\Omega.$
\end{lemma}

\section{Regional strategic actuators}
\label{rbcsection3}
\setcounter{section}{3}\setcounter{equation}{0}
The characteristic of actuators to achieve the regionally approximately boundary controllable of
the system $(\ref{problem})$ will be  explored in this section.

As cited in \cite{AEIJAI},  a actuator can be expressed by a couple $(D,g)$ where
$D\subseteq \Omega $ is the support of the actuator and $g$ is its spatial distribution.
To state our main results,  it is supposed that the control are made by $p$ actuators $(D_i,g_i)_{1\leq i\leq p}$ and
let $Bu=\sum\limits_{i=1}^p{p_{D_i}g_i(x)u_i(t)}$, where $p\in \mathbf{N}$, $g_i(x)\in Z$, $u=(u_1,u_2,\cdots,u_p)$ and $u_i(t)\in L^2(0,b)$.
As cited in \cite{chen2}, all these
distributed parameter systems with moving sensors and actuators form the so-called cyber-physical systems,
 which  are rich in real world applications.  For instance,
in the pest spreading process, $p$ is the number the spreading machines and $u_i(\cdot)$ stands for the control
input strategic of every spreading machines with respect to time $t$ \cite{caljx}.
Then the system $(\ref{problem})$ can be rewritten as
\begin{equation}\label{problem3}
\left\{\begin{array}{l}
_0D^{\alpha}_{t}z(t,x)=Az(t,x)+\sum\limits_{i=1}^p{p_{D_i}g_i(x)u_i(t)},~  (t,x)\in [0,b]\times \Omega,
\\ \lim\limits_{t\to 0^+} {}_0I^{1-\alpha}_{t}z(t,x)=z_0(x).
\end{array}\right.
\end{equation}
Moreover, we suppose that $-A$ is a self-adjoint uniformly elliptic operator,
by \cite{Hilbert}, we get that
there exists a sequence $(\lambda_j,\xi_{jk}): k=1,2,\cdots,r_j$, $j=1,2,\cdots $
such that

$(1)$ For each $ j=1,2,\cdots$, $\lambda_j$ is the eigenvalue of operator $A$ with multiplicities $r_j$ and
\[0>\lambda_1>\lambda_2> \cdots> \lambda_j>\cdots, ~~\lim\limits_{j\to \infty}\lambda_j=-\infty.\]

$(2)$ For each $ j=1,2,\cdots$,  $\xi_{jk}(k=1,2,\cdots,r_j)$ is the orthonormal eigenfunction corresponding to $\lambda_j$, i.e.,
\[
(\xi_{jk_m},\xi_{jk_n})=\left\{\begin{array}{l}1,~~k_m=k_n,
\\
0,~~k_m\neq k_n,
\end{array}\right.\]
where $1\leq k_m,k_n\leq  r_j,$ $k_m,k_n\in \mathbf{N}$ and  $(\cdot,\cdot)$ is the inner product of space $Z$.

Hence, the sequence $\{\xi_{jk}, k=1,2,\cdots,r_j, j=1,2,\cdots \}$ is
a  orthonormal basis in $Z$, the strongly continuous semigroup $\{\Phi(t)\}_{t\geq 0}$ on  $Z$ generated by $A$ is
\begin{eqnarray}\label{phi}
\Phi(t)z(x)=\sum\limits_{j=1}^{\infty}{\sum\limits_{k=1}^{r_j}{\exp(\lambda_jt)(z,\xi_{jk})\xi_{jk}(x)}}, ~~x\in \Omega
\end{eqnarray}
and for any $z(x)\in Z$, it can be expressed as
$
z(x)=\sum\limits_{j=1}^{\infty}{\sum\limits_{k=1}^{r_j}(z,\xi_{jk})\xi_{jk}(x)}.
$
\begin{definition}
A actuators (suite of actuators) is said to be $\Gamma-$strategic if the system under consideration is
regionally approximately boundary controllable on $\Gamma$  at time $b$.
\end{definition}

Before to show our main result in this part, by  Eq.$(\ref{phi})$, for any $z\in L^2(\Omega),$ we have
\begin{eqnarray*}
K_\alpha(t)z(x)&=&{\alpha}\int_0^\infty{\theta\phi_{\alpha}(\theta)\Phi(t^{\alpha}\theta)z(x)}d\theta\\
&=&{\alpha}\int_0^\infty{\theta\phi_{\alpha}(\theta)\sum\limits_{j=1}^{\infty}{\sum\limits_{k=1}^{r_j}
{\exp(\lambda_jt^{\alpha}\theta)(z,\xi_{jk})\xi_{jk}(x)}}}d\theta\\
&=&\sum\limits_{j=1}^{\infty}{\sum\limits_{k=1}^{r_j}\sum\limits_{n=0}^{\infty}\frac{\alpha(\lambda_jt^{\alpha})^n}{n!}}
(z,\xi_{jk})\xi_{jk}(x)\int_0^\infty{\theta^{n+1}\phi_{\alpha}
}d\theta\\
&=&\sum\limits_{j=1}^{\infty}{\sum\limits_{k=1}^{r_j}\sum\limits_{n=0}^{\infty}\frac{\alpha(n+1)!(\lambda_jt^{\alpha})^n}{\Gamma(\alpha n+\alpha+1)n!}}
(z,\xi_{jk})\xi_{jk}(x)\\
&=&\sum\limits_{j=1}^{\infty}\sum\limits_{k=1}^{r_j}\alpha E_{\alpha,\alpha+1}^2(\lambda_jt^{\alpha})(z,\xi_{jk})\xi_{jk}(x),
\end{eqnarray*}
where $E_{\alpha,\beta}^\mu(z):=\sum\limits_{n=0}^{\infty}\frac{(\mu)_n}{\Gamma(\alpha n+\beta)}\frac{z^n}{n!}$, $z\in \mathbf{C}$,
$\alpha,\beta,\mu \in \mathbf{C}$, $\mathbf{Re}{\kern 2pt}\alpha>0$ is the generalized Mittag-Leffler function in three parameters
and here, $(\mu)_n$ is the Pochhammer symbol
defined by (see \cite{H-function}, Section 2.1.1)
\begin{eqnarray}(\mu)_n=\mu(\mu+1)\cdots(\mu+n-1),~n\in \mathbf{N}.\end{eqnarray}
 If $\alpha,\beta \in \mathbf{C} $  such that $\mathbf{Re}{\kern 2pt}\alpha>0$, $\mathbf{Re}{\kern 2pt} \beta>1$, then
(see Section 2.3.4, \cite{Mathai-Haubold} or Section 5.1.1, \cite{Gorenflo})
\begin{eqnarray}
\alpha E_{\alpha,\beta}^2=E_{\alpha,\beta-1}-(1+\alpha-\beta)E_{\alpha,\beta}.
\end{eqnarray}
It then follows that
\begin{eqnarray}\label{3.5}
K_\alpha(t)z(x)
=\sum\limits_{j=1}^{\infty}\sum\limits_{k=1}^{r_j} E_{\alpha,\alpha}(\lambda_jt^{\alpha})(z,\xi_{jk})\xi_{jk}(x)
\end{eqnarray}
and
\begin{eqnarray}
\int_0^t\tau^{\alpha-1}{K_\alpha(\tau)Bu(t-\tau)}d\tau
=\sum\limits_{j=1}^\infty\sum\limits_{k=1}^{r_j}
\sum\limits_{i=1}^p\int_0^t{g^i_{jk}u_i(t-\tau)\tau^{\alpha-1}E_{\alpha,\alpha}
(\lambda_j\tau^\alpha)}d\tau\xi_{jk}(x),
\end{eqnarray}
where $g^i_{jk}=(p_{D_i}g_i,\xi_{jk})$, $j=1,2,\cdots $, $k=1,2,\cdots,r_j$, $i=1,2,\cdots,p$ and
$E_{\alpha,\beta}(z):=\sum\limits_{i=0}^\infty
{\frac{z^i}{\Gamma(\alpha i+\beta)}},$ $\mathbf{Re}{\kern 2pt}\alpha>0, ~\beta,z\in \mathbf{C}$
 is known as the generalized Mittag-Leffler function in two parameters.

\begin{theorem} For any $j=1,2,\cdots$,
 define $p\times r_j$ matrices $G_j$ as
\begin{equation} \label{G}
G_j=\left[ {\begin{array}{*{20}{c}}
{g_{j1}^1}&{g_{j2}^1}&{\cdots}&{g_{jr_j}^1}\\
{g_{j1}^2}&{g_{j2}^2}&{\cdots}&{g_{jr_j}^2}\\
{\vdots}&{\vdots}&{\vdots}&{\vdots}\\
{g_{j1}^p}&{g_{j2}^p}&{\cdots}&{g_{jr_j}^p}
\end{array}} \right],\end{equation}
where $g^i_{jk}=(p_{D_i}g_i,\xi_{jk})$, $j=1,2,\cdots $, $k=1,2,\cdots,r_j$, $i=1,2,\cdots,p$.
 Then the suite of actuators $(D_i,g_i)_{1\leq i\leq p}$ is said to be $\Gamma-$strategic  if and only if
\begin{eqnarray}
p\geq r=\max\{r_j\}~\mbox{ and } ~rank ~G_j=r_j \mbox{ for } j=1,2,\cdots.
\end{eqnarray}
\end{theorem}
$\mathbf{Proof.}$
For any $z_*\in H^{\frac{1}{2}}(\Gamma)$, denote by $\left(\cdot,\cdot\right)_{H^{1/2}(\Gamma)}$ the inner product of space $H^{\frac{1}{2}}(\Gamma),$  we then see that
\begin{eqnarray}\label{attainaletransform}
\left(p_\Gamma \gamma Hu, z_*\right)_{H^{1/2}(\Gamma)}=\sum\limits_{j=1}^\infty\sum\limits_{k=1}^{r_j}
\sum\limits_{i=1}^p\int_0^b{\tau^{\alpha-1}E_{\alpha,\alpha}
(\lambda_j\tau^\alpha)u_i(b-\tau)}d\tau g^i_{jk}z_{jk}=0,~t\in [0,b],
\end{eqnarray}
where $z_{jk}=\left(p_\Gamma \gamma\xi_{jk},z_*\right)_{H^{1/2}(\Gamma)}$, $j=1,2,\cdots$, $k=1,2,\cdots,r_j$.
Further,  Lemma $\ref{infnitydifferential}$  gives
\begin{eqnarray}\label{attainaletransform2}
\sum\limits_{j=1}^\infty\sum\limits_{k=1}^{r_j}t^{\alpha-1}E_{\alpha,\alpha}
(\lambda_jt^\alpha) g^i_{jk}z_{jk}=\textbf{0}_p:=(0,0,\cdots,0)\in \mathbf{R}^p \mbox{ for } t>0,  i=1,2,\cdots, p.
\end{eqnarray}
Then we  conclude that the suite of actuators $(D_i,g_i)_{1\leq i\leq p}$ is $\Gamma-$strategic if and only if
 \begin{eqnarray}\label{3.10}
\sum\limits_{j=1}^{\infty}{
b^{\alpha-1}E_{\alpha,\alpha}
(\lambda_jb^\alpha)}G_jz_j=\textbf{0}_p\Rightarrow z_*=0,
\end{eqnarray}
where  $z_j=(z_{j1},z_{j2},\cdots,z_{jr_j})^T$
is a vector in $\mathbf{R}^{r_j}$ and $j=1,2,\cdots$.

$(a)$ If we assume that  $p\geq r=\max\{r_j\}$ and $rank ~G_j<r_j \mbox{ for some } j=1,2,\cdots$,  there exists a nonzero
element $\tilde{z} \in H^{\frac{1}{2}}(\Gamma)$ with
$\tilde{z}_{j}=\left(\tilde{z}_{j1},\tilde{z}_{j2},\cdots, \tilde{z}_{jr_j}\right)^T\in \mathbf{R}^{r_j}$ such that
\begin{eqnarray}
G_{j}\tilde{z}_{j}=\mathbf{0}_p.
\end{eqnarray}
It then follows from  $E_{\alpha,\alpha}
(\lambda_jt^\alpha)>0 ~(t\geq 0)$   that we can find a nonzero vector  $\tilde{z} $ satisfying
 \begin{eqnarray}
\sum\limits_{j=1}^{\infty}{
b^{\alpha-1}E_{\alpha,\alpha}
(\lambda_jb^\alpha)}G_j\tilde{z}_j=\textbf{0}_p.
\end{eqnarray}
This means that the actuators $(D_i,f_i)_{1\leq i\leq p}$ are not $\Gamma-$strategic.

$(b)$  However, on the contrary, if the actuators $(D_i,g_i)_{1\leq i\leq p}$ are not $\Gamma-$strategic, i.e.,  $\overline{Im(p_\Gamma \gamma  H )}\neq H^{\frac{1}{2}}(\Gamma),$
then there exists a nonzero element  $z\neq \mathbf{0}_n$ satisfying
\begin{eqnarray}
\left(p_\Gamma \gamma  Hu,z\right)_{H^{1/2}(\Gamma)}=0 \mbox{ for all } u\in L^2(0,b; \mathbf{R}^p).
\end{eqnarray}
Then we can find a nonzero element $z_{j^*}\in \mathbf{R}^{r_j}$ such that
\begin{eqnarray}
G_{j^*}z_{j^*}=\mathbf{0}_p.
\end{eqnarray}
This allows us to complete the conclusion of the theorem.

\section{Regional boundary controllability with minimum energy control}
\label{rbcsection4}
\setcounter{section}{4}\setcounter{equation}{0}

In this section, we explore the possibility of finding a minimum energy control when
the system $(\ref{problem})$ can be steered from a given initial vector $z_0$ to a
target function $z_b$ on the boundary subregion $\Gamma.$ The method used here
is an extension of those in \cite{AEIJAI,AEIJAI2,217,zerrik,179}.

Consider the following minimization problem
\begin{eqnarray}\label{minimum}
\left\{\begin{array}{l}
\inf\limits_u J(u)=\int_0^b{\|u(t)\|^2_{\mathbf{R}^p}}dt
\\ u\in U_b=\{u\in  L^2\left(0,b; \mathbf{R}^p\right): p_\Gamma \gamma z(b,u)=z_b\},
\end{array}\right.
\end{eqnarray}
where, obviously, $U_b$ is a closed convex set.
We then show a direct approach to the solution of the minimum energy problem $(\ref{minimum})$.

\begin{theorem}\label{theorem3.2}
If the system $(\ref{problem})$ is  regionally approximately boundary controllable on $\Gamma$, then for any
$z_b\in H^{\frac{1}{2}}(\Gamma),$ the minimum energy problem $(\ref{minimum})$ has a unique solution given by
\begin{eqnarray}u^*(t)=\left(p_\Gamma \gamma H\right)^*R_\Gamma^{-1} \left(z_b-p_\Gamma \gamma K_\alpha(b)z_0\right),\end{eqnarray}
where $R_\Gamma=p_\Gamma \gamma HH^*\gamma ^*p^*_\Gamma$ and $H^*$ is defined in Eq.$(\ref{H*})$.
\end{theorem}
$\textbf{Proof.}$ To begin with, since the solution of $(\ref{problem})$ excited by the control $u^*$ is given by
\begin{eqnarray}
z(t,u^*)=K_\alpha(t)z_0+\int_0^t(t-s)^{\alpha-1}{K_\alpha(t-s)Bu^*(s)}ds,
\end{eqnarray}
we get that
\begin{eqnarray*}
p_\Gamma \gamma z(b,u^*)&=&p_\Gamma \gamma\left[ K_\alpha(b)z_0+\int_0^b(b-s)^{\alpha-1}{K_\alpha(b-s)Bu^*(s)}ds\right]\\
&=&p_\Gamma \gamma K_\alpha(b)z_0+p_\Gamma \gamma H \left(p_\Gamma \gamma H\right)^*R_\Gamma^{-1} \left(z_b-p_\Gamma \gamma K_\alpha(b)z_0\right)
\\&=&z_b.
\end{eqnarray*}

Next, we show that if the system $(\ref{problem})$ is  regionally approximately boundary
 controllable on $\Gamma$  at time $b$, then the operator $R_\Gamma$ is coercive. In fact, for any $z_1\in H^{\frac{1}{2}}(\Gamma)$,
 there exists a control $u\in L^2(0,b,\mathbf{R}^p)$ such that
 \begin{eqnarray}
z_1=p_\Gamma \gamma \left[K_\alpha(b)z_0+Hu\right]
\end{eqnarray}
 and
\begin{eqnarray*}\left<R_\Gamma z_1,z_1\right>_{H^{1/2}(\Gamma)}&=&
\left\|H^*\gamma ^*p^*_\Gamma z_1\right\|^2_{L^2(0,b,\mathbf{R}^p)} \\&=&
\left\|B^*(b-\cdot)^{\alpha-1}K_\alpha^*(b-\cdot)\gamma ^*p^*_\Gamma z_1\right\|^2_{L^2(0,b,\mathbf{R}^p)}\\&\geq&
\left\| z_1\right\|^2_{H^{1/2}(\Gamma)}
.\end{eqnarray*}
Moreover, since $R_\Gamma\in \mathcal{L}\left(H^{\frac{1}{2}}(\Gamma),H^{\frac{1}{2}}(\Gamma)\right)$,
by the Theorem 1.1 in \cite{Lions2},  it follows that $R_\Gamma$ is an isomorphism.

Finally, we prove that  $u^*$ solves the minimum energy problem $(\ref{minimum})$.
For this purpose, since $p_\Gamma \gamma  z(b,u^*)=z_b,$ for any $u\in L^2(0,b,\mathbf{R}^p )$
with $p_\Gamma \gamma z(b,u)=z_b$,
one has
\begin{eqnarray}p_\Gamma \gamma \left[z(b,u^*)-z(b,u)\right]=0,\end{eqnarray}
which follows that
\begin{eqnarray*}
0&=&p_\Gamma \gamma \int_0^b{(b-s)^{\alpha-1}K_\alpha(b-s)B \left[u^*(s)-u(s)\right]}ds
=p_\Gamma \gamma H\left[u^*-u\right].\end{eqnarray*}
 Thus, by
\begin{eqnarray*}J'(u^*)(u^*-u)&=&
2\int_0^b{\left<u^*(s)-u(s),u^*(s)\right>}ds\\
&=&2\int_0^b{\left<u^*(s)-u(s),\left(p_\Gamma \gamma H\right)^*R_\Gamma^{-1} \left(z_b-p_\Gamma \gamma K_\alpha(b)z_0\right)\right>}ds\\
&=&2\int_0^b{\left<p_\Gamma \gamma H\left[u^*(s)-u(s)\right],R_\Gamma^{-1} \left(z_b-p_\Gamma \gamma K_\alpha(b)z_0\right)\right>}ds\\&=&0,\end{eqnarray*}
 it  follows that $J(u)\geq J(u^*)$, i.e., $u^*$ solves the minimum energy problem $(\ref{minimum})$ and the proof is complete.
\section{The connection between  internal and boundary regional controllability}
\label{rbcsection5}
\setcounter{section}{5}\setcounter{equation}{0}

Based on an intension of the regional controllability of integer order differential equations
developed in \cite{217,zerrik},  we here give a transfer on the internal and boundary
regional controllability of  fractional order sub-diffusion equations $(\ref{problem})$ and
 develop two types of controls, i.e., zone or pointwise.

\subsection{Internal and boundary regional controllability}
In this part, we present a internal and boundary
regional controllability transfer of  the problem $(\ref{problem})$.
To this end, suppose that $z(b,u)\in Z$ and we first define a operator
\begin{eqnarray} T: H^{\frac{1}{2}}(\partial \Omega) \to H^1(\Omega) \mbox{ such that } \gamma Tg=g.~\forall g\in H^{\frac{1}{2}}(\partial \Omega),\end{eqnarray}
which is linear and continuous \cite{Lions}. Let $z_b\in H^{\frac{1}{2}}(\Gamma)$ with the extension
$p_\Gamma^* z_b \in H^{\frac{1}{2}}(\partial \Omega)$
and consider the sets
\begin{eqnarray}
\Omega_1=\left\{Tp_\Gamma^* z_b\in  Z| z_b\in H^{\frac{1}{2}}(\Gamma) \right\}
\mbox{ and } \Omega_2=\mathop  \cup \limits_{z_b\in H^{1/2}(\Gamma)}\mbox{Supp }Tp_\Gamma^* z_b.
\end{eqnarray}
For any $r>0$ be arbitrary sufficiently small, consider
\begin{eqnarray}
D_r=\mathop  \cup \limits_{z\in\Gamma} B(z,r) \mbox{ and let  }  \omega_r=D_r \cap \Omega_2,
\end{eqnarray}
 where $B(z,r)$ is a ball of radius $r$ centred in $z$.
\begin{theorem}\label{theorem5.1}
If the system $(\ref{problem})$ is exactly(respectively, approximately) controllable on $\omega_r$,
then it is also  exactly(respectively, approximately) boundary controllable on $\Gamma.$
\end{theorem}
\textbf{ Proof.}
Let $z_b\in H^{\frac{1}{2}}(\Gamma)$ be the target function. By utilizing the trace theorem
\cite{Retherford}, there exists $Tp_\Gamma^* z_b\in  Z$ with a bounded support such that
$\gamma (Tp_\Gamma^* z_b)=p_\Gamma^* z_b.$ Then

$1)$ if the system $(\ref{problem})$ is exactly controllable on $\omega_r$, for any $y_b\in H^1(\omega_r)$,
there exists a control $u\in L^2(0,b;\mathbf{R}^p)$ such that
\begin{eqnarray}
p_{\omega_r} z(b,u)=y_b.
\end{eqnarray}
Then $p_{\omega_r} T p_\Gamma^* z_b\in H^1(\omega_r)$ and
there exists a control $u\in L^2(0,b;\mathbf{R}^p)$ such that
\begin{eqnarray}
p_{\omega_r} z(b,u)=p_{\omega_r} T p_\Gamma^* z_b\mbox{ and } \gamma p_{\omega_r} z(b,u)=p_\Gamma^* z_b.
\end{eqnarray}
Thus
$
p_\Gamma \gamma p_{\omega_r} z(b,u)=z_b,
$
i.e., the system $(\ref{problem})$ is exactly boundary controllable on $\Gamma.$

$2)$ if the system $(\ref{problem})$ is approximately controllable on $\omega_r$, for and $\varepsilon>0$  and any $y_b\in H^1(\omega_r)$,
there exists a control $u\in L^2(0,b;\mathbf{R}^p)$ such that
\begin{eqnarray}
\left\|p_{\omega_r} z(b,u)-y_b\right\|_{H^1(\omega_r)}\leq \varepsilon.
\end{eqnarray}
Then  for any $\varepsilon>0,$
there exists a control $u\in L^2(0,b;\mathbf{R}^p)$ such that
\begin{eqnarray}
\left\|p_{\omega_r} z(b,u)-p_{\omega_r} T p_\Gamma^* z_b\right\|_{H^1(\omega_r)} \leq \varepsilon.
\end{eqnarray}
Moreover, by the continuity of the trace mapping $\gamma$ on $H^1(\omega_r)$, one has
\begin{eqnarray}
\left\|\gamma(p_{\omega_r} z(b,u))-\gamma(p_{\omega_r} T p_\Gamma^* z_b)\right\|_{H^1(\partial\omega_r)} \leq \varepsilon,
\end{eqnarray}
therefore
$
\left\|p_\Gamma \gamma(p_{\omega_r} z(b,u))- z_b\right\|_{H^1(\Gamma)} \leq \varepsilon,
$
Thus  $(\ref{problem})$ is approximately boundary controllable on $\Gamma$ and the proof is complete.

\subsection{Regional boundary target control}

This part is concerned with the approach for the control which drives the problem $(\ref{problem})$ from $z_0$ to $z_b$ on $\Gamma.$
Let $z_b\in H^{\frac{1}{2}}(\Gamma)$ with the extension
$p_\Gamma^* z_b \in H^{\frac{1}{2}}(\partial \Omega)$. By Theorem $\ref{theorem5.1}$, the problem may be solved by driving the
system $(\ref{problem})$ from $z_0$ to $y_b\in H^1(\omega_r)$ on $\omega_r.$

The following two sets will be used in our discussion.
\begin{eqnarray}
G=\{g\in H^1(\Omega): g=0\mbox{ in }\Omega\backslash\omega_r \} \mbox{ and } E=\{e\in H^1(\Omega):e =0\mbox{ in } \omega_r \}.
\end{eqnarray}

\subsubsection{Case of zone actuator}
Let us consider the system $(\ref{problem})$ with a zone actuator $(D,f)$ where $D\subseteq \Omega$ is
the support of the actuator and $f$ is its spatial distribution. Then the system can be written in the
form
\begin{eqnarray}\label{problemzone}
\left\{\begin{array}{l}
_0D^{\alpha}_{t}z(x,t)=Az(x,t)+p_Df(x)u(t)~\mbox{ in }~\Omega\times [0,b],
\\  \lim\limits_{t\to 0^+} {}_0I^{1-\alpha}_{t}z(x,t)=z_0(x)~\mbox{ in }~\Omega,
\\z(x, t)=0~\mbox{ on }~\partial\Omega\times [0,b].
\end{array}\right.
\end{eqnarray}
For any $g\in G$, consider the system
\begin{equation}\label{5.13}
\left\{\begin{array}{l}
Q{}_tD^{\alpha}_b\left[(b-t)^{1-\alpha}\varphi(x,t)\right]=A^*Q\left[(b-t)^{1-\alpha}\varphi(x,t)\right]~\mbox{ in }~\Omega\times [0,b],
\\ \lim\limits_{t\to 0^+} Q{}_tI^{1-\alpha}_b\left[(b-t)^{1-\alpha}\varphi(x,t)\right]=p_{\omega_r}^*g(x)~\mbox{ in }~\Omega,
\\ \varphi (x, t)=0~\mbox{ on }~\partial\Omega\times [0,b].
\end{array}\right.
\end{equation}
where $Q$ is a reflection operator on interval $[0, b]$ such that
\begin{eqnarray}
Qf(t) := f( b - t).
\end{eqnarray}
By the argument in \cite{Klimek}, we see that the following properties on operator $Q$ hold:
\begin{eqnarray}
Q_0I_t^\alpha f(t)={}_tI_b^\alpha Qf(t),~~~Q{}_0D_t^\alpha f(t)={}_tD_b^\alpha Qf(t)\end{eqnarray}
and
\begin{eqnarray}\label{5.14}
_0I_t^\alpha Qf(t) = Q{}_tI_b^\alpha f(t),~~~{}_0D_t^\alpha Qf(t)=Q{}_tD_b^\alpha f(t) .\end{eqnarray}
Then system \textcolor[rgb]{1.00,0.00,0.00}{$(\ref{5.13})$} can be rewritten as
\begin{equation}
\left\{\begin{array}{l}
{}_0D^{\alpha}_tQ\left[(b-t)^{1-\alpha}\varphi(x,t)\right]=A^*Q\left[(b-t)^{1-\alpha}\varphi(x,t)\right]~\mbox{ in }~\Omega\times [0,b],
\\ \lim\limits_{t\to 0^+} {}_0I^{1-\alpha}_tQ\left[(b-t)^{1-\alpha}\varphi(x,t)\right]=p_{\omega_r}^*g(x)~\mbox{ in }~\Omega,
\\ \varphi (x, t)=0~\mbox{ on }~\partial\Omega\times [0,b]
\end{array}\right.
\end{equation}
and its unique mild solution is
$\varphi(x,t)=(b-t)^{\alpha-1}K_\alpha^*(b-t)p_{\omega_r}^*g(x).$
Moreover, we define the semi-norm
\begin{equation}\label{Gnorm}
g\in G \to \|g\|_G^2=\int_0^b{ \left(f,\varphi (\cdot,t)\right)_{L^2(D)}^2}dt
\end{equation}
on $G$ and obtain the following result.
\begin{lemma}\label{norm}
$(\ref{Gnorm})$ defines a norm on $G$ if the system $(\ref{problemzone})$ is  regionally approximately controllable on $\omega$ at time $b$.
\end{lemma}
\textbf{Proof.}
For any $g\in G$, if the system $(\ref{problemzone})$ is regionally approximately controllable on $\omega$,  we have
\begin{eqnarray*}
Ker(  H^* p^*_\omega) =Ker\left[(b-s)^{\alpha-1}\left(f,K_\alpha^*(b-s)p_{\omega_r}^* g\right)_{L^2(D)}\right]
=Ker\left[\left(f,\varphi (\cdot,t)\right)_{L^2(D)}\right]=\{0\}.
\end{eqnarray*}
It then follows from
  \begin{eqnarray*}
\|g\|_G^2=\int_0^b{ \left(f,\varphi (\cdot,t)\right)_{L^2(D)}^2}dt=0
\Leftrightarrow \left(f,\varphi (\cdot,t)\right)_{L^2(D)}=0\end{eqnarray*}
 that    $\|\cdot\|_G$ is a norm of space $G$ and the proof is complete.

Moreover, let $u(t)=\left(f,\varphi (\cdot,t)\right)_{L^2(D)}$ and decomposed the system
$(\ref{problemzone})$ into an autonomous system and a homogeneous initial condition one
\begin{equation}
\left\{\begin{array}{l}
_0D^{\alpha}_{t}\psi_1(x,t)=A\psi_1(x,t)+p_Df(x)\left(f,\varphi (\cdot,t)\right)_{L^2(D)}~\mbox{ in }~ \Omega\times [0,b],
\\ \lim\limits_{t\to 0^+} {}_0I^{1-\alpha}_t\psi_1(x,t)=0~~\mbox{ in }~\Omega,
\\ \psi_1 (x, t)=0~\mbox{ on }~\partial\Omega\times [0,b]
\end{array}\right.
\end{equation}
and
\begin{eqnarray}
 \left\{\begin{array}{l}
_0D^{\alpha}_{t}\psi_2(x,t)=A\psi_2(x,t)~\mbox{ in }~ \Omega\times [0,b],
\\ \lim\limits_{t\to 0^+} {}_0I^{1-\alpha}_t\psi_2(x,t)=z_0(x)~\mbox{ in }~ \Omega,
\\ \psi_2 (x, t)=0~\mbox{ on }~\partial\Omega\times [0,b].
\end{array}\right.
\end{eqnarray}
Let $\wedge$ be the operator $\wedge: G\to E^\bot$ given by
\begin{eqnarray}\label{5.19}
\wedge g=p_{\omega_r} \psi_1(\cdot,b), ~~\forall g\in G.
\end{eqnarray}
Then for any $z_b\in H^1(\omega_r),$   the regional control problem on $\omega_r$ is equivalent to the resolution of the equation
\begin{eqnarray}\label{4.13}
\wedge g= z_b- p_{\omega_r}\psi_2(\cdot,b)
\end{eqnarray}
and we have the following result.

\begin{theorem}\label{Th2}
Assume that the system $(\ref{problemzone})$ is  regionally approximately controllable on $\omega_r$ at time $b$,
 then $(\ref{4.13})$ admits a unique solution $g\in G$ and the control
\begin{eqnarray}u^*(t)=\left(f,\varphi (\cdot,t)\right)_{L^2(D)}\end{eqnarray}
steers the problem $(\ref{problemzone})$ to $z_b$ on $\omega_r$. Moreover, $u^*$ solves the minimum energy problem
 \begin{eqnarray}\label{costf}
\inf\limits_u J(u)=\int_0^b{\|u(t)\|^2_{\mathbf{R}^p}}dt.
\end{eqnarray}
\end{theorem}
\textbf{Proof.}
From Lemma $\ref{norm}$, if the system $(\ref{problemzone})$ is
regionally approximately controllable on $\omega_r$ at time $b$, then $\|\cdot\|_G$ is a norm of space $G$.
 Let the completion of $G$ with respect to the norm $\|\cdot\|_G$ again by $G$.

Next,  we show that $(\ref{4.13})$ admits a unique solution in $G$.
For any $g\in G$, by Eq. $(\ref{5.19})$, it follows that
\begin{eqnarray*}\left<g,\wedge g\right>&=&\left<g, p_{\omega_r} \psi_1(\cdot,b)\right>\\&=&
\left<g,p_{\omega_r}\int_0^b(b-s)^{\alpha-1}{K_\alpha(b-s)p_Df(\cdot)\left(f,\varphi (\cdot,s)\right)_{L^2(D)}}ds\right>\\&=&
\int_0^b{\| \left(f,\varphi (\cdot,t)\right)_{L^2(D)}\|^2}ds=
\|g\|^2_G.\end{eqnarray*}
Hence,  it follows from the Theorem 1.1 in \cite{Lions2}
that $(\ref{4.13})$ admits a unique solution in $G$.

Let $u=u^*$ in problem $(\ref{problemzone})$, then $p_{\omega_r} z(b,u^*)=z_b.$ Finally, we show that $u^*$ minimize
the const functional $(\ref{costf}).$
For any $u_1\in L^2(0,b,\mathbf{R}^p )$ with $p_{\omega_r} z(b,u_1)=z_b$, we have
\begin{eqnarray}p_{\omega_r}\left[z(b,u^*)-z(b,u_1)\right]=0.\end{eqnarray}
Then
\begin{eqnarray*}
0=p_{\omega_r}\int_0^b{(b-s)^{\alpha-1}K_\alpha(b-s)p_Df(x) \left[u^*(s)-u_1(s)\right]}ds.\end{eqnarray*}
Moreover, since
\begin{eqnarray*}J'(u^*)(u^*-u_1)&=&
2\int_0^b{(u^*(s)-u_1(s))u^*(s)}ds\\
&=&2\int_0^b{(u^*(s)-u_1(s))\left(f,\varphi (\cdot,t)\right)_{L^2(D)}}ds\\
&=&2\int_0^b{\left(p_Df\left[u^*(s)-u_1(s)\right],(b-t)^{\alpha-1}K_\alpha^*(b-t)p_{\omega_r}^*g\right)}ds \\
&=&2\left(p_{\omega_r}\int_0^b{(b-s)^{\alpha-1}K_\alpha(b-s)p_Df(x) \left[u^*(s)-u_1(s)\right]}ds,g\right)\\&=&0,\end{eqnarray*}
one has $J(u)\geq J(u^*)$, i.e.,
$u^*$ solves the minimum energy problem $(\ref{costf})$ and the proof is complete.

\subsubsection{Case of pointwise actuator}
Consider the system $(\ref{problem})$ with a pointwise internal  actuator, which can be written in the form
\begin{eqnarray}\label{problempointwise}
\left\{\begin{array}{l}
_0D^{\alpha}_{t}z(x,t)=Az(x,t)+\delta(x-\sigma)u(t)~\mbox{ in }~\Omega\times [0,b],
\\  \lim\limits_{t\to 0^+} {}_0I^{1-\alpha}_{t}z(x,t)=z_0~\mbox{ in }~\Omega,
\\z(x, t)=0~\mbox{ on }~\partial\Omega\times [0,b],
\end{array}\right.
\end{eqnarray}
where $\sigma$  is the actuator support.
For any $g\in G$, consider \textcolor[rgb]{1.00,0.00,0.00}{  $(\ref{5.13})$ }and
define the semi-norm
\begin{eqnarray}g\to \|g\|^2_G=\int_0^b{\left\|\varphi(\sigma,s)\right\|^2}ds,\end{eqnarray} which defines a norm on $G$
if $(\ref{problempointwise})$ is regionally approximately controllable.

Similar to the argument in section 5.2.1, let $u(t)=\varphi (\sigma,t)$ and we consider the following system
\begin{equation}
\left\{\begin{array}{l}
_0D^{\alpha}_{t}\psi_1(x,t)=A\psi_1(x,t)+\delta(x-\sigma)\varphi (\sigma,t)~\mbox{ in }~  \Omega\times [0,b],
\\ \lim\limits_{t\to 0^+} {}_0I^{1-\alpha}_t\psi_1(x,t)=0~\mbox{ in }~ \Omega,
\\ \psi_1 (x, t)=0~\mbox{ on }~ \partial\Omega\times [0,b]
\end{array}\right.
\end{equation}
and
\begin{eqnarray}
 \left\{\begin{array}{l}
_0D^{\alpha}_{t}\psi_2(x,t)=A\psi_2(x,t)~\mbox{ in }~\Omega\times [0,b],
\\ \lim\limits_{t\to 0^+} {}_0I^{1-\alpha}_t\psi_2(x,t)=z_0(x)~\mbox{ in }~\Omega,
\\ \psi_2 (x, t)=0~\mbox{ on }~\partial\Omega\times [0,b].
\end{array}\right.
\end{eqnarray}
Then the regional control problem on $\omega_r$ is equivalent to the resolution of the equation
\begin{eqnarray}\label{5.28}
\wedge g= z_b- p_{\omega_r}\psi_2(\cdot,b)
\end{eqnarray}
and we see the following result.

\begin{theorem}\label{Th2}
Assume that the system $(\ref{problempointwise})$ is  regionally approximately controllable on $\omega_r$ at time $b$,
 then $(\ref{5.28})$ admits a unique solution $g\in G$ and the control
\begin{eqnarray}\label{pointwisecontrolinut}
u^*(t)=\varphi (\sigma,t)
\end{eqnarray}
steers  $(\ref{problemzone})$ to $z_b$ on $\omega_r$. Moreover, this control minimize the cost functional $(\ref{costf})$.
\end{theorem}
\subsubsection{Simulation}
The resolution of the regional boundary control problem may be seen via the
following simplified steps (see the case of pointwise actuator for example).

1) Initial data $\Omega$, $\Gamma,$ $z_b$ and the actuator;

2) Solve the problem $(\ref{5.28})$ $(\to$ $g)$;

3) Solve the problem $(\ref{5.13})$ $(\to$ $\varphi(\sigma,t))$;

4) Apply the control $u^*(t) = \varphi(\sigma,t)$.

For example, consider the system $(\ref{examplel})$ and
let $\Omega=[0,1]\times [0,1],$ $\Gamma=\{0\}\times [1/4,3/4]$, $b=5$.  For the target function
$z_b$ on $\Gamma \subseteq \partial \Omega$,  we assume that
\begin{eqnarray}
z_b(0,y) =\left\{\begin{array}{l}
0,~~~~~~0\leq y<1/4;
\\ 0.017+4(y-1/4)^2(y-3/4)^2,~~~~1/4\leq y\leq 3/4;
\\0,~~~~~~3/4<y\leq 1
\end{array}\right.
\end{eqnarray}
and the actuator is supposed to be located in
$D=\{0 \}\times \{0.5\}\subseteq \Omega$.

\begin{figure}
\begin{minipage}[t]{1\linewidth}
\centering
\includegraphics[width=0.8\textwidth]{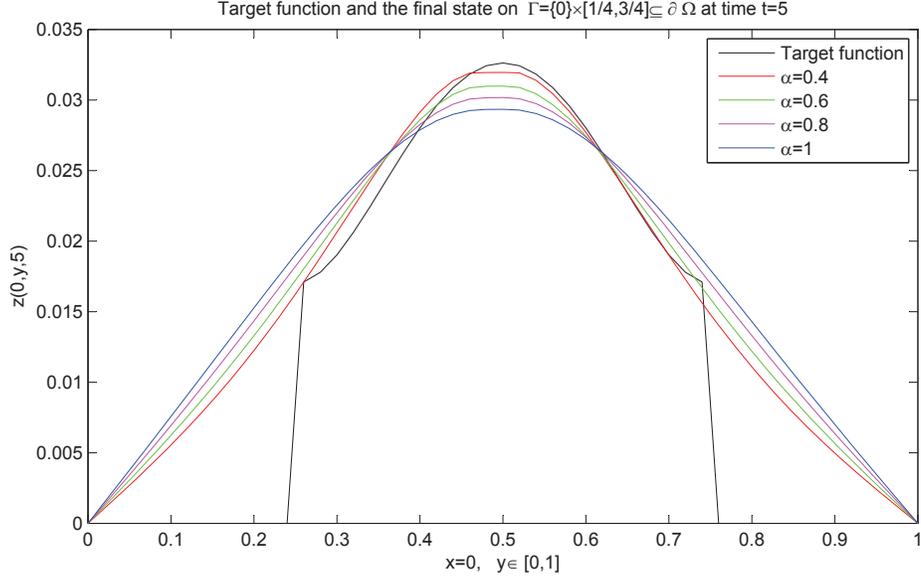}
\end{minipage}%
\caption{Final reached state and  target function
on $\Gamma \subseteq \partial \Omega$ at time $t=5$.}
\label{fig:ft}
\end{figure}

\begin{figure}
\begin{minipage}[t]{1\linewidth}
\centering
\includegraphics[width=0.8\textwidth]{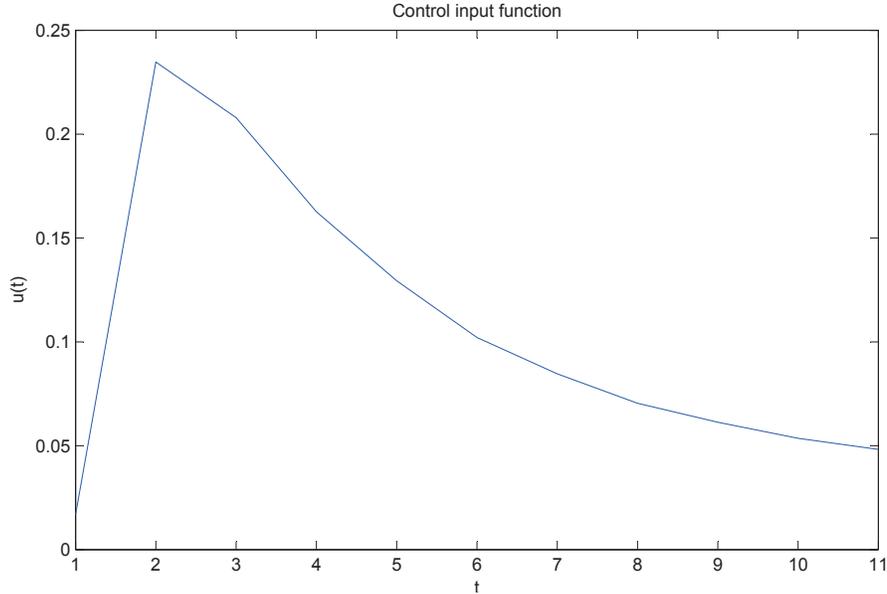}
\end{minipage}%
\caption{Control input function, which is calculated by the formula $(\ref{pointwisecontrolinut}).$}
\label{fig:controlinput}
\end{figure}

Figure $\ref{fig:ft}$ shows how the final reached state is very close to the target function
on $\Gamma \subseteq \partial \Omega$ at time $t=5$ when $\alpha=0.4,~0.6,~0.8,~1.0$. This also implies that
time fractional  diffusion systems can offer better performance compared with those using integer order
distributed parameter systems. Moreover,  when $\alpha=0.4$, the corresponding control input, which is
calculated by the formula $(\ref{pointwisecontrolinut}),$ is presented at Figure $\ref{fig:controlinput}.$ 


\section{CONCLUSIONS}
In this paper,  the regional boundary controllability of the Riemann-Liouville time fractional diffusion systems
of order $\alpha\in (0,1]$ is discussed, which is motivated by many realistic situation encountered in various applications.
The results here provide some insights into the qualitative analysis of the design of fractional
order diffusion equations, which can also be extended to complex fractional order distributed
parameter dynamic systems. Various open questions are still under consideration. The problem of
constrained control as well as the case of fractional order distributed parameter dynamic systems
with more complicated regional sensing and actuation configurations are of great interest.
For more
information on the potential topics related to fractional order distributed parameter systems, we refer the readers to \cite{ge2015JAS}
and the references therein.
\\

\noindent {\bf Acknowledgement.}

This work was supported  by Chinese Universities Scientific Fund (No.CUSF-DH-D-2014061)
 and the Natural Science Foundation of Shanghai (No. 15ZR1400800).

\bibliographystyle{model1b-num-names}

\begin{thebibliography}{00}
\addtolength\itemsep{-0.38cm}

\bibitem{AEIJAI} El Jai A. and Pritchard A. J., {\it Sensors and controls in the analysis of distributed systems.}
 Halsted Press, New York, 1988.

\bibitem{AEIJAI2} El Jai A., Distributed systems analysis via sensors and actuators. {\it Internat. J. of Sensors
and Actuators} 29 (1991), 1-11.


\bibitem{217} Zerrik E., Boutoulout A. and El Jai A., Acutators and regional boundary controllability for parabolic systems.
{\it Internat. J. Systems Sci.} 31 (2000), 73-82.

\bibitem{zerrik} Zerrik E., Badraoui L. and  El Jai A., Sensors and regional boundary state reconstruction of parabolic systems.
{\it Sensors and Actuators } 75 (1999), 102-117.

\bibitem{179} Sakawa Y.,  Controllability for partial differential equations of parabolic type. {\it Siam J. Control } 12 (1974), 389-400.

\bibitem{CTRWfirst} Montroll E. W. and  Weiss G. H., Random walks on lattices. II. {\it J. Mathematical Phys.} 6 (1965), 167-181.


\bibitem{ralf} Metzler R. and  Klafter J., The random walk's guide to anomalous diffusion: a fractional dynamics approach.
{\it Phys. Rep.} 339 (2000), 1-77.

\bibitem{ralf2} Metzler R. and Klafter J.,
The restaurant at the end of the random walk: recent developments in the description of anomalous transport by fractional dynamics.
{\it J. Phys. A } 37 (2004), 161-208.

 \bibitem{Cartea} Cartea $\acute{A}$. and Negrete D., Fluid limit of the continuous-time random walk
 with general L$\acute{e}$vy jump distribution functions. {\it Phys. Rev. E } 76  041105 (2007), 1-8.

\bibitem{3M} Meerschaert M. M. and Scalas E., Coupled continuous time random walks in finance. {\it Phys. A } 370 (2006),  114-118.

\bibitem{Adams} Adams E. E. and  Gelhar L.W., Field study of dispersion in a heterogeneous aquifer 2. Spatial moments analysis.
{\it Water Resources Res. } 28 (1992), 3293-3307.

\bibitem{Uchaikin} Uchaikin V. and Sibatov R.,
{\it Fractional Kinetics in Solids: Anomalous Charge Transport in Semiconductors.}
 Published by World Scientific Publishing Co. Pte. Ltd., 2013.

\bibitem{2Spears} Spears W. M. and Spears D. F., {\it Physicomimetics: Physics-based swarm intelligence.}
Springer Science $\&$ Business Media, 2012.


 \bibitem{opAgrawal} Agrawal  O. P.,
Solution for a fractional diffusion-wave equation defined in a bounded domain. {\it Nonlinear Dynam. } 29 (2002), 145-155.

\bibitem{Podlubny} Podlubny I., {\it Fractional differential equations.}  Academic Press, San Diego, 1999.

\bibitem{Kilbas} Kilbas A. A., Srivastava H. M. and Trujillo J. J., {\it Theory and applications of fractional differential equations.}
 Elsevier, 2006.

\bibitem{Torvik}
Torvik P.~J.  and Bagley R.~L., On the appearance of the fractional derivative in
  the behavior of real materials, {\it J. Appl. Mech.} 51~(2) (1984), 294--298.

\bibitem{Mandelbrot}
Mandelbrot B.~B., {\it  The fractal geometry of nature.} Vol. 173, Macmillan, 1983.

\bibitem{GeFDTA-rc} Ge F., Chen Y. and Kou C., Regional controllability of anomalous diffusion generated by
the time fractional diffusion equations, {\it In: ASME IDETC/CIE 2015, Boston, Aug. 2-5,
2015, DETC2015-46697. See also: arXiv preprint arXiv:1508.00047.}

\bibitem{Mainardi} Mainardi F., Paradisi P. and Gorenflo R., {\it Probability distributions generated by fractional diffusion equations.}
 in: J. Kertesz, I. Kondor (Eds.), Econophysics: An Emerging Science, Kluwer, Dordrecht, 2000.

\bibitem{zhouyong} Zhou Y. and Jiao F., Existence of mild solutions for fractional neutral evolution equations.
{\it Comput. Math. Appl.} 59 (2010),  1063-1077.

 \bibitem{dualityrelationship} Pritchard A. J. and Wirth A.,
Unbounded control and observation systems and their duality.
 {\it Siam J. Control Optim.} 16 (1978), 535-545.

\bibitem{Bernard} Dacorogna  B., {\it Direct methods in the calculus of variations.
 Second edition. } Appl. Math. Sci. 78. Springer, New York, 2008.

\bibitem{chen2}  Coopmans C., Stark B., Jensen A., Chen Y. and McKee M.,  {\it Cyber-physical systems enabled
by small unmanned aerial vehicles: a chapter in Handbook of Unmanned Aerial
Vehicles.} Valavanis, Kimon P.; Vachtsevanos, George J (Eds.), 2015, Springer.

\bibitem{caljx} Cao J., Chen Y. and Li C., Multi-UAV-based optimal crop-dusting of anomalously diffusing
infestation of crops. {\it arXiv preprint arXiv:1411.2880 (2014).}

 \bibitem{Hilbert} Courant  R. and Hilbert D., {\it Methods of mathematical physics. } Vol. I. Interscience, New York, 1966.

\bibitem{H-function} Erd$\acute{e}$lyi, A., Magnus, W., Oberhettinger, F. and Tricomi, F. G., {\it Higher transcendental functions, Vol. I.}
McGraw-Hill, NewYork-Toronto-London, 1953.

\bibitem{Mathai-Haubold} Mathai A. M. and  Haubold H. J.,
{\it Special functions for applied scientists. } Springer Science Business Media, LLC.  2008.

\bibitem{Gorenflo} Gorenflo R, Kilbas A. A.,
Mainardi F and  Rogosin S. V., {\it Mittag-Leffler functions, related topics and applications.} Springer-Verlag Berlin Heidelberg, 2014.

\bibitem{Lions2}
 Lions J.~L., {\it Optimal control of systems governed by partial differential
  equations}, Vol. 170, Springer Verlag, 1971.

\bibitem{Lions} Lions J. L., Exact controllability, stabilization and perturbations for distributed systems.
{\it Siam Rev.}  30 (1988),  1-68.


\bibitem{Klimek} Ma${\l}$gorzata K., {\it On solutions of linear fractional differential equations of a variational type.}
Publishing Office of Czestochowa University of Technology. 2009.

\bibitem{Retherford}  Retherford J. R., {\it Hilbert space: compact operators and the trace theorem.}
London Mathematical Society Student Texts, 27. Cambridge University Press, Cambridge, 1993.

\bibitem{ge2015JAS}  F. Ge, Y. Chen, C. Kou, Cyber-physical systems as general distributed parameter
systems: three types of fractional order models and emerging research opportunities
[J]. {\it IEEE/CAA J. Autom. Sin.} 2 (4) (2015), 353-357.

\end{thebibliography}

\end{document}